\documentclass[12pt]{amsart}
\usepackage{amsfonts}
\usepackage{amsmath,amssymb,amsthm}

\newcommand{\R}{\mathbb R}

\newcommand{\C}{\mathbb C}

\renewcommand{\span}{\mathrm{span}}
\newcommand{\tr}{\mathrm{tr}}

\newtheorem{thm}{Theorem}[section]
\newtheorem{cor}[thm]{Corollary}

\newtheorem{prop}[thm]{Proposition}
\theoremstyle{definition}
\newtheorem{defn}[thm]{Definition}
\theoremstyle{remark}

\newcommand{\ds}{\displaystyle}

\begin{document}

\title[INVARIANTS AND BONNET-TYPE THEOREM FOR SURFACES IN $\R^4$]
{INVARIANTS AND BONNET-TYPE THEOREM FOR SURFACES IN $\R^4$}

\author{Georgi Ganchev and Velichka Milousheva}
\address{Bulgarian Academy of Sciences, Institute of Mathematics and Informatics,
Acad. G. Bonchev Str. bl. 8, 1113 Sofia, Bulgaria}
\email{ganchev@math.bas.bg}

\address{Bulgarian Academy of Sciences, Institute of Mathematics and Informatics,
Acad. G. Bonchev Str. bl. 8, 1113, Sofia, Bulgaria; "L. Karavelov"
Civil Engineering Higher School, 175 Suhodolska Str., 1373 Sofia,
Bulgaria} \email{vmil@math.bas.bg}



\begin{abstract}
In the tangent plane at any point of a surface in the
four-dimensional Euclidean space we consider an invariant linear
map of Weingarten-type and find a geometrically determined moving
frame field. Writing derivative formulas of Frenet-type for this
frame field, we obtain eight invariant functions. We prove a
fundamental theorem of Bonnet-type, stating that these eight
invariants under some natural conditions determine the surface up
to a motion.

We show that the basic geometric classes of surfaces in the
four-dimensional Euclidean space, determined by conditions on
their invariants, can be interpreted in terms of the properties of
two geometric figures: the tangent indicatrix, which is a conic in
the tangent plane, and the normal curvature ellipse.

We construct a family of surfaces with flat normal connection.
\end{abstract}

\maketitle

\section{Introduction}

Local invariants of surfaces in the four-dimensional Euclidean
space $\R^4$ were studied by Eisenhart \cite{Ein},  Kommerell
\cite{Kom}, Moore and Wilson \cite{GM2}, Schouten and Struik
\cite{Scho-Str}, Spivak \cite{Spivak}, Wong \cite{Wong}, Little
\cite{Little}. Their study was based on a special configuration,
namely a point and an ellipse lying in the normal space (the ellipse of normal curvature).
This configuration leads to a
theory of axial principal directions, along which the vector-valued second
fundamental form points in the direction of the major and the
minor axes of the curvature ellipse. In higher dimensions there is
also a similar configuration consisting of a point and a Veronese
manifold. This configuration determines second order scalar
invariants and generates principal axes ''in general''
\cite{Little}. Points where the construction of principal axes fails are regarded as
singularities of the field of axes. Geometric
singularities  for immersions in Riemannian manifolds
are considered in \cite{Asp}. Special types of tangent vector fields on a
surface in $\R^4$ are defined in terms of the properties of the
normal curvature ellipse and families of lines determined by such
tangent vector fields are studied in \cite{Gar-Soto,Mello}.

In this paper our aim is to develop the local theory of surfaces
in $\R^4$ on the base of the Weingarten map similarly to the
classical case of surfaces in $\R^3$.

Let $M^2$ be a surface in $\R^4$ with tangent space $T_pM^2$ at
any point $p \in M^2$. In \cite{GM1} we introduced an invariant
linear map $\gamma$ of Weingarten-type at any $T_pM^2$, which
generates two invariant functions $k$ and $\varkappa$. The sign of
the function $k$ is a geometric invariant and the sign of
$\varkappa$ is invariant under motions in $\R^4$. However, the
sign of $\varkappa$ changes under symmetries with respect to a
hyperplane in $\R^4$. Analogously to $\R^3$, the invariants $k$
and $\varkappa$ divide the points of $M^2$ into four types: flat,
elliptic, hyperbolic and parabolic points. The surfaces consisting
of flat points, i.e. satisfying the condition $k= \varkappa = 0$,
either lie in $\R^3$ or are developable ruled surfaces in $\R^4$.
Everywhere in the present considerations, we exclude the points at
which $k= \varkappa = 0$.

The minimal surfaces in $\R^4$ are characterized in terms of the
invariants $k$ and $\varkappa$ by the condition $\varkappa^2 - k
=0$, and the surfaces with flat normal connection are
characterized by $\varkappa = 0$.

Further, the map $\gamma$ generates the corresponding second
fundamental form $II$ at any point $p \in M^2$ in the standard
way. In \cite{GM3} we gave a geometric interpretation of the
second fundamental form $II$ of the surface. We introduced an
invariant $\zeta_{g_1,\,g_2}$ of a pair of two tangents $g_1$,
$g_2$ at any point $p$ of $M^2$. Then the tangents $g_1$, $g_2$
are conjugate in terms of $II$ if and only if
$\zeta_{g_1,\,g_2}=0$. The notions of a normal curvature and a
geodesic torsion of a tangent were introduced by means of the
invariant $\zeta$. It turns out that asymptotic tangents and
principal tangents in terms of $II$  are characterized by zero
normal curvature and zero geodesic torsion, respectively. The
principal normal curvatures $\nu'$ and $\nu''$ are defined as the
normal curvatures of the principal tangents. The invariants $k$
and $\varkappa$ satisfy the equalities
$$k = \nu' \nu''; \qquad \varkappa = \frac{\nu' + \nu''}{2}.$$

It turns out that the points at which any tangent is principal
(''umbilical" points) are characterized by zero mean curvature
vector, i.e. the surfaces consisting of "umbilical" points are
exactly the minimal surfaces in $\R^4$.

The indicatrix of Dupin at an arbitrary (non-flat) point of a
surface in $\R^3$ is introduced by means of the second fundamental
form. Here, we introduce the indicatrix $\chi$ at any point $p \in
M^2$ in the same way:
$$\chi : \nu' X^2 + \nu'' Y^2 = \varepsilon, \qquad \varepsilon = \pm1.$$
The conjugacy in terms of the second fundamental form coincides
with the conjugacy with respect to the indicatrix $\chi$.

In Section 3 we prove that: \vskip 1mm \emph{The surface $M^2$ is
minimal if and only if the indicatrix $\chi$ is a circle.} \vskip
1mm \emph{The surface $M^2$ is with flat normal connection if and
only if the indicatrix $\chi$ is a rectangular hyperbola (a
Lorentz circle).} \vskip 1mm

In the local theory of surfaces a statement of significant
importance is a theorem of Bonnet-type giving the natural
conditions under which the surface is determined up to a motion. A
theorem of this type was proved for surfaces with flat normal
connection by B.-Y. Chen in \cite{Chen1}. In the class of surfaces
with $\varkappa^2 - k >0$ we find a geometrically determined
moving frame of Frenet-type. Considering the corresponding
derivative formulas, we obtain eight invariant functions. In
Section 4 we prove our basic Theorem \ref{T:Main Theorem}, stating
that \vskip 1mm \emph{The eight invariant functions, satisfying
some natural conditions, determine the surface up to a motion in
$\R^4$.} \vskip 1mm

In Section 5 we construct a family of surfaces with flat normal
connection, which lie on a standard rotational hypersurface in
$\R^4$, and describe those of them with constant Gauss curvature,
constant mean curvature, and constant invariant $k$.

\section{A geometric interpretation of the second fundamental form}

Let $M^2: z = z(u,v), \, \, (u,v) \in {\mathcal D}$ (${\mathcal D}
\subset \R^2$) be a 2-dimensional surface in $\R^4$. The tangent
space $T_pM^2$ to $M^2$ at an arbitrary point $p=z(u,v)$ of $M^2$
is ${\rm span} \{z_u, z_v\}$.   We choose an orthonormal normal
frame field $\{e_1, e_2\}$ of $M^2$ so that the quadruple $\{z_u,
z_v, e_1, e_2\}$ is positive oriented in $\R^4$. Then the
following derivative formulas hold:
$$\begin{array}{l}
\vspace{2mm} \nabla'_{z_u}z_u=z_{uu} = \Gamma_{11}^1 \, z_u +
\Gamma_{11}^2 \, z_v
+ c_{11}^1\, e_1 + c_{11}^2\, e_2,\\
\vspace{2mm} \nabla'_{z_u}z_v=z_{uv} = \Gamma_{12}^1 \, z_u +
\Gamma_{12}^2 \, z_v
+ c_{12}^1\, e_1 + c_{12}^2\, e_2,\\
\vspace{2mm} \nabla'_{z_v}z_v=z_{vv} = \Gamma_{22}^1 \, z_u +
\Gamma_{22}^2 \, z_v
+ c_{22}^1\, e_1 + c_{22}^2\, e_2,\\
\end{array}$$
where $\Gamma_{ij}^k$ are the Christoffel's symbols and
$c_{ij}^k$, $i, j, k = 1,2$ are functions on $M^2$.

We use the standard denotations \;$E, \; F, \; G$ for the
coefficients of the first fundamental form and set
$W=\sqrt{EG-F^2}$. Denoting by $\sigma$ the second fundamental
tensor of $M^2$, we have
$$\begin{array}{l}
\sigma(z_u,z_u)=c_{11}^1\, e_1 + c_{11}^2\, e_2,\\
[2mm]
\sigma(z_u,z_v)=c_{12}^1\, e_1 + c_{12}^2\, e_2,\\
[2mm] \sigma(z_v,z_v)=c_{22}^1\, e_1 + c_{22}^2\,
e_2.\end{array}$$ The three pairs of normal vectors
$\{\sigma(z_u,z_u), \sigma(z_u,z_v)\}$, $\{\sigma(z_u,z_u),
\sigma(z_v,z_v)\}$, $\{\sigma(z_u,z_v), \sigma(z_v,z_v)\}$ form
three parallelograms with oriented areas
$$\Delta_1 = \left|%
\begin{array}{cc}
\vspace{2mm}
  c_{11}^1 & c_{12}^1 \\
  c_{11}^2 & c_{12}^2 \\
\end{array}%
\right|, \quad
\Delta_2 = \left|%
\begin{array}{cc}
\vspace{2mm}
  c_{11}^1 & c_{22}^1 \\
  c_{11}^2 & c_{22}^2 \\
\end{array}%
\right|, \quad
\Delta_3 = \left|%
\begin{array}{cc}
\vspace{2mm}
  c_{12}^1 & c_{22}^1 \\
  c_{12}^2 & c_{22}^2 \\
\end{array}%
\right|,$$ respectively. These oriented areas determine three
functions $\ds{L = \frac{2 \Delta_1}{W}, \,\, M = \frac{
\Delta_2}{W},\,\, N = \frac{2 \Delta_3}{W}}$, which change in the
same way as the coefficients $E, F, G$ under any change of the
parameters $(u,v)$.

Using the functions $E$, $F$, $G$ and $L$, $M$, $N$, in \cite{GM1}
we introduced the linear map $\gamma$ in the tangent space at any
point of $M^2$
$$\gamma: T_pM^2 \rightarrow T_pM^2,$$
defined by the equalities
$$\begin{array}{l}
\vspace{2mm}
\gamma(z_u)=\gamma_1^1z_u+\gamma_1^2z_v,\\
\vspace{2mm} \gamma(z_v)=\gamma_2^1z_u+\gamma_2^2z_v,
\end{array}$$
where
$$\displaystyle{\gamma_1^1=\frac{FM-GL}{EG-F^2}, \quad
\gamma_1^2 =\frac{FL-EM}{EG-F^2}}, \quad
\displaystyle{\gamma_2^1=\frac{FN-GM}{EG-F^2}, \quad
\gamma_2^2=\frac{FM-EN}{EG-F^2}}.$$

The linear map $\gamma$ of Weingarten type at the point $p \in
M^2$ is invariant with respect to changes of parameters on $M^2$
as well as to motions in $\R^4$. This implies that the functions
$$k = \frac{LN - M^2}{EG - F^2}, \qquad
\varkappa =\frac{EN+GL-2FM}{2(EG-F^2)}$$ are invariants of the
surface $M^2$.

The invariant $\varkappa$ turns out to be the curvature of the
normal connection of the surface $M^2$ in $\R^4$.

As in the classical case, the invariants $k$ and $\varkappa$
divide the points of $M^2$ into four types: flat, elliptic,
parabolic and hyperbolic. The surfaces consisting of flat points
satisfy the conditions
$$k(u,v)=0, \quad \varkappa(u,v)=0, \qquad (u,v) \in \mathcal D,$$
or equivalently $L(u,v)=0, \; M(u,v)=0, \; N(u,v)=0, \; (u,v) \in
\mathcal D.$ These surfaces are either planar surfaces (there
exists a hyperplane $\R^3 \subset \R^4$ containing $M^2$) or
developable ruled surfaces in $\R^4$.

Further we consider surfaces free of flat points, i.e. $(L, M, N)
\neq (0, 0, 0)$.

\vskip 2mm Let $X = \lambda z_u + \mu z_v, \,\, (\lambda,\mu) \neq
(0,0)$ be a tangent vector at a point $p \in M^2$. The Weingarten
map $\gamma$ determines a second fundamental form of the surface
$M^2$ at $p$ as follows \cite{GM1}:
$$II(\lambda,\mu) = - g(\gamma(X),X) = L\lambda^2 + 2M\lambda\mu + N\mu^2, \quad \lambda,\mu \in \R.$$

As in the classical differential geometry of surfaces in $\R^3$
the second fundamental form $II$ determines conjugate tangents at
a point $p$ of $M^2$. Two tangents $g_1: X = \lambda_1 z_u + \mu_1
z_v$ and $g_2: X = \lambda_2 z_u + \mu_2 z_v$ are said to be
\textit{conjugate tangents}  if $II(\lambda_1, \mu_1; \lambda_2,
\mu_2) = 0$, i.e.
$$L\lambda_1 \lambda_2 + M (\lambda_1 \mu_2 +\lambda_2 \mu_1) + N\mu_1 \mu_2 = 0.$$

A geometric interpretation of the second fundamental form and the
map $\gamma$ can be given using the geometric approach in
\cite{GM3} to the notion of conjugacy.

Conjugate tangents are introduced in \cite{GM3} in a geometric way
as follows. Let $g$ be a tangent at the point $p \in M^2$
determined by the vector $X =  \lambda z_u + \mu z_v$. We consider
the linear map $\sigma_g: T_pM^2 \rightarrow (T_pM^2)^{\bot}$,
defined by
$$\sigma_g(Y) = \ds{\sigma \left(\frac{\lambda z_u + \mu z_v}{\sqrt{I(\lambda, \mu)}},\, Y \right)}, \quad Y \in T_pM^2. $$

Let $g_1: X_1 =  \lambda_1 z_u + \mu_1 z_v$ and $g_2: X_2 =
\lambda_2 z_u + \mu_2 z_v$ be two tangents at the point $p \in
M^2$. The oriented areas of the parallelograms spanned by the
pairs of normal vectors $\sigma_{g_1} (z_u)$, $\sigma_{g_2} (z_v)$
and $\sigma_{g_2} (z_u)$, $\sigma_{g_1} (z_v)$  are denoted by
$S(\sigma_{g_1} (z_u), \sigma_{g_2} (z_v))$, and $S(\sigma_{g_2}
(z_u), \sigma_{g_1} (z_v))$, respectively. Then
 $$\zeta_{\,g_1,g_2} = \ds{\frac{S(\sigma_{g_1} (z_u),\sigma_{g_2} (z_v))}{W} + \frac{S(\sigma_{g_2} (z_u), \sigma_{g_1} (z_v))}{W}}$$
is an invariant of the tangents $g_1$, $g_2$.

\begin{defn}\cite{GM3} \label{D:conjugate}
Two tangents $g_1: X_1 = \lambda_1 z_u + \mu_1 z_v$ and $g_2: X_2
= \lambda_2 z_u + \mu_2 z_v$ are said to be \emph{conjugate
tangents}   if $\zeta_{\,g_1,g_2} = 0$.
\end{defn}

Calculating the oriented areas in $\zeta_{\,g_1,g_2}$, it can be
found that
$$\zeta_{\,g_1,g_2} = \ds{\frac{L \lambda_1 \lambda_2 + M (\lambda_1 \mu_2 + \mu_1 \lambda_2) + N \mu_1 \mu_2}{\sqrt{I(\lambda_1, \mu_1)} \sqrt{I(\lambda_2, \mu_2)}}}
= \ds{\frac{II(\lambda_1, \mu_1; \lambda_2,
\mu_2)}{\sqrt{I(\lambda_1, \mu_1)} \sqrt{I(\lambda_2, \mu_2)}}}.$$
Thus, $\zeta_{\,g_1,g_2} = 0$ if and only if $L\lambda_1 \lambda_2
+ M (\lambda_1 \mu_2 +\lambda_2 \mu_1) + N\mu_1 \mu_2 = 0.$ Hence,
the tangents $g_1$ and $g_2$ are conjugate according to Definition
\ref{D:conjugate} if and only if they are conjugate with respect
to  the second fundamental form $II$.

We defined two invariants $\nu_g$ and $\alpha_g$ of any tangent
$g$ of the surface in terms of $\zeta_{\,g_1,g_2}$ as follows:
$$\nu_g = \zeta_{\,g,g}; \qquad \alpha_g = \zeta_{\,g, g^{\bot}}.$$
We call $\nu_g$ the \emph{normal curvature} of  $g$, and
$\alpha_g$ - the \emph{geodesic torsion} of $g$. The invariant
$\nu_g$ is expressed by the first and the second fundamental forms
of the surface in the same way as the normal curvature of a
tangent in the theory of surfaces in $\R^3$, i.e. $\nu_g =
\ds{\frac{II(\lambda, \mu)}{I(\lambda, \mu)}}$. Further, the
invariant $\alpha_g$ can be written in the following way:
$$\alpha_g = \ds{\frac{\lambda^2 (EM - FL) + \lambda \mu (EN - GL) + \mu^2(FN - GM)}{W I(\lambda, \mu)}}.$$
Hence, $\alpha_g$ is expressed by the coefficients of the first
and the second fundamental forms in the same way as the geodesic
torsion in the theory of surfaces in $\R^3$.

The notions of asymptotic tangents and principal tangents are
defined in terms of the conjugacy as in $\R^3$.

A tangent $g: X = \lambda z_u + \mu z_v$ is said to be
\textit{asymptotic}  if it is self-conjugate, i.e. $L\lambda^2 +
2M\lambda\mu + N\mu^2 = 0$. Hence, a tangent $g$ is asymptotic if
and only if $\nu_g = 0$. If $p$ is an elliptic point of $M^2$
($k>0$) then there are no asymptotic tangents through $p$; if $p$
is a hyperbolic point ($k<0$)  then there are two asymptotic
tangents passing through $p$, and if $p$ is a parabolic point
($k=0$) then there is one asymptotic tangent through $p$. Thus,
the sign of the invariant $k$ determines the number of asymptotic
tangents at the point.

A tangent $g: X = \lambda z_u + \mu z_v$ is said to be
\textit{principal} if it is perpendicular to its conjugate. The
equation for the principal tangents at a point $p \in M^2$ is
$$\left|\begin{array}{cc}
E & F\\
[2mm] L & M \end{array}\right| \lambda^2+ \left|\begin{array}{cc}
E & G\\
[2mm] L & N \end{array}\right| \lambda \mu+
\left|\begin{array}{cc}
F & G\\
[2mm] M & N \end{array}\right| \mu^2=0.$$ Hence, a tangent $g$ is
principal if and only if $\alpha_g = 0$. As in the classical case
we have $\varkappa^2 - k \geq 0$ at each point of the surface. If
$\varkappa^2 - k =0$,  every tangent is principal, and if
$\varkappa^2 - k >0$, there exist exactly two principal tangents.

A line $c: u=u(q), \; v=v(q); \; q\in \textrm{J} \subset \R$ on
$M^2$ is said to be an \textit{asymptotic line}  if its tangent at
any point is asymptotic. A line $c: u=u(q), \; v=v(q); \; q\in
\textrm{J} \subset \R$ on $M^2$ is said to be a \textit{principal
line} (a \textit{line of curvature}) if its tangent at any point
is principal. The surface $M^2$ is parameterized by principal
lines if and only if $F=0,\,\, M=0$.

\vskip 2mm Let $M^2$ be a surface with $\varkappa^2 - k>0$ at each
point. We assume that $M^2$ is parameterized by principal lines
and denote the unit vector fields
$\displaystyle{x=\frac{z_u}{\sqrt E}}$, $\ds{y=\frac{z_v}{\sqrt
G}}$. The equality $M = 0$ implies that the normal vector fields
$\sigma(x,x)$ and $\sigma(y,y)$ are collinear. We denote by $b$ a
unit normal vector field collinear with $\sigma(x,x)$ and
$\sigma(y,y)$, and by $l$ the unit normal vector field such that
$\{x,y,b,l\}$ is a positive oriented orthonormal frame field of
$M^2$ (the two vectors $\{b, l\}$ are determined up to a sign).
Thus we obtain a geometrically determined orthonormal frame field
$\{x,y,b,l\}$ at each point $p \in M^2$. With respect to the frame
field $\{x,y,b,l\}$ we have the following formulas:
$$\begin{array}{l}
\vspace{2mm}
\sigma(x,x) = \nu_1\,b; \\
\vspace{2mm}
\sigma(x,y) = \lambda\,b + \mu\,l;  \\
\vspace{2mm} \sigma(y,y) = \nu_2\,b,
\end{array}\leqno{(2.1)}$$
where $\nu_1, \nu_2, \lambda, \mu$ are invariant functions, whose
signs depend on the orientation of $\{b, l\}$.

Hence the invariants $k$, $\varkappa$, and the Gauss curvature $K$
of $M^2$ are expressed as follows:
$$k = - 4\nu_1\,\nu_2\,\mu^2, \quad \quad \varkappa = (\nu_1-\nu_2)\mu, \quad \quad
K = \nu_1\,\nu_2- (\lambda^2 + \mu^2).\leqno(2.2)$$ Since
$\varkappa^2 - k > 0$, equalities (2.2) imply that $\mu \neq 0$.

The normal mean curvature vector field  of $M^2$ is $H =
\ds{\frac{\sigma(x,x) + \sigma(y,y)}{2} = \frac{\nu_1 +
\nu_2}{2}\, b}$. Taking into account (2.2) we obtain that the
length $\Vert H \Vert$ of the mean curvature vector field is given
by the formula
$$\Vert H \Vert = \displaystyle{\frac{\sqrt{\varkappa^2-k}}{2 |\mu |}},$$
which shows that $|\mu|$ is expressed by the invariants $k$,
$\varkappa$ and the mean curvature function.

Now we shall discuss the geometric meaning of the invariant
$\lambda$. Let $M$ be an $n$-dimensional submanifold of
$(n+m)$-dimensional Riemannian manifold $\widetilde{M}$ and $\xi$
be a normal vector field of $M$. In \cite{Chen1} B.-Y. Chen
defined the \emph{allied vector field} $a(\xi)$ of $\xi$ by the
formula
$$a(\xi) = \ds{\frac{\|\xi\|}{n} \sum_{k=2}^m \{\tr(A_1 A_k)\}\xi_k},$$
where $\{\xi_1 = \ds{\frac{\xi}{\|\xi\|}},\xi_2, \xi_m \}$ is an
orthonormal base of the normal space of $M$, and $A_i = A_{\xi_i},
\,\, i = 1,\dots, m$ is the shape operator with respect to
$\xi_i$. In particular, the allied vector field $a(H)$ of the mean
curvature vector field $H$ is a well-defined normal vector field
which is orthogonal to $H$. It is called the \emph{allied mean
curvature vector field} of $M$ in $\widetilde{M}$. B.-Y. Chen
defined  the $\mathcal{A}$-submanifolds to be those submanifolds
of $\widetilde{M}$ for which
 $a(H)$ vanishes identically \cite{Chen1}.
In \cite{GVV1,GVV2} the $\mathcal{A}$-submanifolds are called
\emph{Chen submanifolds}. It is easy to see that minimal
submanifolds, pseudo-umbilical submanifolds and hypersurfaces are
Chen submanifolds. These Chen submanifolds are said to be trivial
$\mathcal{A}$-submanifolds. Now let $M^2$ be a surface in $\R^4$.
Applying the definition of the allied mean curvature vector field
from equalities (2.1) we get
$$a(H) = \ds{\frac{\nu_1 + \nu_2}{2} \,\lambda \mu \, l}= \ds{\frac{\sqrt{\varkappa^2-k}}{2} \, \lambda \,l}.$$
Hence, if $M^2$ is free of minimal points, then $a(H) = 0$ if and
only if $\lambda = 0$. This gives the geometric meaning of the
invariant $\lambda$. It is clear that $M^2$ is a non-trivial  Chen
surface if and only if the invariant $\lambda$ is zero.

 \vskip 2mm
The normal curvatures $\nu' = \ds{\frac{L}{E}}$ and $\nu'' =
\ds{\frac{N}{G}}$ of the principal tangents are said to be
\textit{principal normal curvatures} of $M^2$. The invariants $k$
and $\varkappa$ of $M^2$ are expressed by the principal normal
curvatures $\nu'$ and $\nu''$ as follows:
$$k = \nu' \nu''; \qquad \varkappa = \frac{\nu' + \nu''}{2}. \leqno{(2.3)}$$

Similarly to the theory of surfaces in $\R^3$, we consider the
indicatrix $\chi$ in the tangent space $T_pM^2$ at an arbitrary
point $p$ of $M^2$, defined by
$$\chi : \nu' X^2 + \nu'' Y^2 = \varepsilon, \qquad \varepsilon = \pm 1.$$

If $p$ is an elliptic point ($k > 0$), then the indicatrix $\chi$
is an ellipse. The axes of $\chi$ are collinear with the principal
directions at the point $p$, and the lengths of the axes are
$\ds{\frac{2}{\sqrt{|\nu'|}}}$ and
$\ds{\frac{2}{\sqrt{|\nu''|}}}$.

If $p$ is a  hyperbolic point ($k < 0$), then the indicatrix
$\chi$ consists of two hyperbolas. For the sake of simplicity we
say that $\chi$ is a hyperbola. The axes of $\chi$ are collinear
with the principal directions, and the lengths of the axes are
$\ds{\frac{2}{\sqrt{|\nu'|}}}$ and
$\ds{\frac{2}{\sqrt{|\nu''|}}}$.

If $p$ is a parabolic point ($k = 0$), then the indicatrix $\chi$
consists of two straight lines parallel to the principal direction
with non-zero normal curvature.

The following statement holds :

\begin{prop}\label{P:conjugacy}
Two tangents $g_1$ and $g_2$ are conjugate tangents of $M^2$ if
and only if $g_1$ and $g_2$ are conjugate with respect to the
indicatrix $\chi$.
\end{prop}

\vskip 2mm
\section{Classes of surfaces characterized in terms of the tangent indicatrix and the normal curvature ellipse}

The minimal surfaces in $\R^4$ are characterized by

\begin{prop}\label{P:minimal} \cite{GM1}
Let $M^2$ be a surface in $\R^4$ free of flat points. Then $M^2$
is minimal if and only if
$$\varkappa^2 - k = 0.$$
\end{prop}

The surfaces with flat normal connection are characterized by

\begin{prop}\label{P:Flat normal}
Let $M^2$ be a surface in $\R^4$ free of flat points. Then $M^2$
is a surface with flat normal connection if and only if
$$\varkappa = 0.$$
\end{prop}

\noindent \emph{Proof:} Let $D$ be the normal connection of $M^2$.
For any tangent vector fields $x,y$ and any normal vector field
$n$ we have the standard decomposition
$$\nabla'_xn = -A_n(x) + D_xn,$$
where $\langle A_n(x), y\rangle = \langle \sigma(x,y),n\rangle$.

The curvature tensor $R^{\bot}$ of the normal connection $D$ is
given by
$$R^{\bot}(x,y)n = D_xD_yn - D_yD_xn - D_{[x,y]}n.$$
Then the normal curvature (the curvature of the normal connection)
at a point $p \in M^2$ is defined by $\langle R^{\bot}(x,y)n_2,
n_1\rangle,$ where $\{x,y,n_1,n_2\}$ is a right oriented
orthonormal quadruple.

Without loss of generality we assume that $F = 0$ and denote the
unit vector fields $\displaystyle{x=\frac{z_u}{\sqrt E}, \;
y=\frac{z_v}{\sqrt G}}$. Then we have the formulas
$$\begin{array}{l}
\vspace{2mm}
\sigma(x,x)=\displaystyle{\frac{c_{11}^1}{E}\;\;n_1 \;+\;\frac{c_{11}^2}{E}\;\;n_2},\\
\vspace{2mm}
\sigma(x,y)= \displaystyle{\frac{c_{12}^1}{\sqrt{EG}}\,n_1 +\frac{c_{12}^2}{\sqrt{EG}}\,n_2,}\\
\vspace{2mm} \sigma(y,y) =\displaystyle{\frac{c_{22}^1}{G}\;\;n_1
\;+\;\;\frac{c_{22}^2}{G}\;\;n_2}.
\end{array} \leqno{(3.1)}$$
Hence,
$$\begin{array}{ll}
\vspace{2mm} A_1(x)=\displaystyle{\frac{c_{11}^1}{E}\;\;x
\;+\;\frac{c_{12}^1}{\sqrt{EG}}\;\;y},
\qquad & A_2(x)=\displaystyle{\frac{c_{11}^2}{E}\;\;x \;+\;\frac{c_{12}^2}{\sqrt{EG}}\;\;y},\\
\vspace{2mm} A_1(y)= \displaystyle{\frac{c_{12}^1}{\sqrt{EG}}\;\;x
+ \frac{c_{22}^1}{G}\;\;y,} \qquad & A_2(y)=
\displaystyle{\frac{c_{12}^2}{\sqrt{EG}}\;\;x +
\frac{c_{22}^2}{G}\;\;y.}
\end{array} \leqno{(3.2)}$$
Using (3.2) we calculate
$$\begin{array}{ll}
\vspace{2mm} (A_2 \circ A_1 - A_1 \circ A_2) (x) =
\ds{\left(\frac{c_{11}^1 c_{12}^2 - c_{11}^2 c_{12}^1}{E
\sqrt{EG}} + \frac{c_{12}^1 c_{22}^2 - c_{12}^2 c_{22}^1}{G
\sqrt{EG}}\right) y} =
\ds{\frac{EN + GL}{2EG}\, y;}\\
\vspace{2mm} (A_2 \circ A_1 - A_1 \circ A_2) (y) = -
\ds{\left(\frac{c_{11}^1 c_{12}^2 - c_{11}^2 c_{12}^1}{E
\sqrt{EG}} + \frac{c_{12}^1 c_{22}^2 - c_{12}^2 c_{22}^1}{G
\sqrt{EG}}\right) x} = - \ds{\frac{EN + GL}{2EG}\, x.}
\end{array}$$
Thus we get
$$\begin{array}{ll}
\vspace{2mm}
(A_2 \circ A_1 - A_1 \circ A_2) (x)= \varkappa \,y;\\
\vspace{2mm} (A_2 \circ A_1 - A_1 \circ A_2) (y)= - \varkappa \,x.
\end{array} \leqno{(3.3)}$$

Note that $A_2 \circ A_1 - A_1 \circ A_2$ is an invariant
skew-symmetric operator in the tangent space, i.e. it does not
depend on the choice of the orthonormal tangent frame field
$\{x,y\}$.

Since the curvature tensor $R'$ of the connection $\nabla'$ is
zero, we have
$$\nabla'_x \nabla'_y n_1 - \nabla'_y \nabla'_x n_1 - \nabla'_{[x,y]} n_1 = 0.$$
Therefore, the tangent component and the normal component of
$R'(x,y)n_1$ are both zero. The normal component is $D_xD_yn_1 -
D_yD_xn_1 - D_{[x,y]}n_1 - \sigma(x,A_1y) + \sigma(y,A_1x).$
Hence,
$$D_xD_yn_1 - D_yD_xn_1 - D_{[x,y]}n_1 = \sigma(x,A_1y) - \sigma(y,A_1x). \leqno{(3.4)}$$
The left-hand side of (3.4) is $R^{\bot}(x,y)n_1$. Then
$$\langle R^{\bot}(x,y)n_1, n_2 \rangle = \langle \sigma(x,A_1y), n_2 \rangle
- \langle \sigma(y,A_1x), n_2 \rangle = \langle (A_2 \circ A_1 -
A_1 \circ A_2) (y), x \rangle.$$ Using (3.3) we obtain $\langle
R^{\bot}(x,y)n_1, n_2 \rangle = - \varkappa$. Since $\langle
R^{\bot}(x,y)n_1, n_2 \rangle = - \langle R^{\bot}(x,y)n_2, n_1
\rangle$, we get
$$\langle R^{\bot}(x,y)n_2, n_1 \rangle = \varkappa.$$
The last equality implies that $\varkappa$ is the curvature of the
normal connection. Hence, $M^2$ is a surface with flat normal
connection if and only if $\varkappa = 0$. \qed

\vskip 4mm We note that the condition $\varkappa = 0$ implies that
$k<0$ and the surface $M^2$ has two families of orthogonal
asymptotic lines.

\vskip 2mm Now we shall characterize the minimal surfaces and the
surfaces with flat normal connection in terms of the tangent
indicatrix of the surface.

\begin{prop}\label{P:Minimal-circle}
Let $M^2$ be a surface in $\R^4$ free of flat points. Then $M^2$
is minimal if and only if at each point of $M^2$ the tangent
indicatrix $\chi$ is a circle.
\end{prop}

\noindent {\it Proof:} Let $M^2$ be a surface in $\R^4$ free of
flat points. From equalities (2.3) it follows that
$$\varkappa^2 - k = \ds{\left(\frac{\nu' - \nu''}{2}\right)^2}.$$
Obviously $\varkappa^2 - k = 0$ if and only if $\nu' = \nu''$.
Applying Proposition \ref{P:minimal}, we get that $M^2$ is minimal
if and only if  $\chi$ is a circle. \qed

\begin{prop}\label{P:Flat normal-Lorentz circle}
Let $M^2$ be a surface in $\R^4$ free of flat points.  Then $M^2$
is a surface with flat normal connection if and only if at each
point of $M^2$ the tangent indicatrix $\chi$ is a rectangular
hyperbola (a Lorentz circle).
\end{prop}

\noindent {\it Proof:} Let $M^2$ be a surface in $\R^4$ free of
flat points. From (2.3) it follows that $\varkappa = 0$ if and
only if $\nu'' = -\nu'$.

If $M^2$ is a surface with flat normal connection, then $k < 0$,
and hence $\chi$ is a hyperbola. From $\nu'' = -\nu'$ it follows
that the semi-axes of $\chi$ are equal to
$\ds{\frac{1}{\sqrt{|\nu'|}}}$, i.e. $\chi$ is a rectangular
hyperbola.

Conversely, if $\chi$ is a rectangular hyperbola, then $\nu'' =
-\nu'$, which implies that $M^2$ is a surface with flat normal
connection. \qed

\vskip 2mm The minimal surfaces and the surfaces with flat normal
connection can also be characterized in terms of the ellipse of
normal curvature.

Let us recall that the \textit{ellipse of normal curvature} at a
point $p$ of a surface $M^2$ in $\R^4$ is the ellipse
$\varepsilon$ in the normal space at the point $p$ given by
$\{\sigma(x,x): \, x \in T_pM^2, \, g(x,x) = 1\}$ \cite{MW1, MW2}.
Let $\{x,y\}$ be an orthonormal base of the tangent space $T_pM^2$
at $p$. Then, for any $v = \cos \psi \, x + \sin \psi \, y$, we
have
$$\sigma(v, v) = H + \ds{\cos 2\psi  \, \frac{\sigma(x,x) - \sigma(y,y)}{2}
+ \sin 2 \psi  \, \sigma(x,y)},$$ where $H = \ds{\frac{\sigma(x,x)
+ \sigma(y,y)}{2}}$ \, is the mean curvature vector of $M^2$ at
$p$. So, when $v$ goes once around the unit tangent circle, the
vector $\sigma(v,v)$ goes twice around the ellipse centered at
$H$. The vectors $\ds{\frac{\sigma(x,x) - \sigma(y,y)}{2}}$ \, and
$\sigma(x,y)$ determine conjugate directions of $\varepsilon$. The
area of the ellipse is given by the formula $S_{\varepsilon} = \pi
\,|S_0|$, where $S_0$ is the oriented area of the parallelogram
spanned by the vectors $\ds{\frac{\sigma(x,x) - \sigma(y,y)}{2}}$
\, and $\sigma(x,y)$. In view of formulas (3.1) the oriented area
$S_0$ is expressed as $S_0 = \ds{\frac{\varkappa}{2}}$. Hence,
$S_{\varepsilon}= \ds{\frac{\pi}{2}\, |\varkappa}|$.

\vskip 2mm A surface $M^2$ in $\R^4$ is called
\textit{super-conformal} \cite{BFLPP}  if at any point of $M^2$
the ellipse of curvature is a circle. In \cite{DT} it is given an
explicit construction of any simply connected super-conformal
surface in $\R^4$ that is free of minimal and flat points. \vskip
2mm Obviously, $M^2$ is minimal if and only if for each point $p
\in M^2$ the ellipse of curvature is centered at $p$.

The  minimal surfaces in $\R^4$ are divided into two subclasses:
\begin{itemize}
\item
the subclass of minimal super-conformal surfaces, characterized by
the condition that the ellipse of curvature is a circle;

\item
the subclass of minimal surfaces of general type, characterized by
the condition that the ellipse of curvature is not a circle.
\end{itemize}

In  \cite{GM2} it is proved that on any minimal surface $M^2$ the
Gauss curvature $K$ and the normal curvature $\varkappa$ satisfy
the following inequality
$$K^2-\varkappa^2\geq 0.$$
The two subclasses of minimal surfaces are characterized in terms
of the invariants $K$ and  $\varkappa$ as follows:
\begin{itemize}
\item the class of minimal super-conformal surfaces is characterized by $K^2 - \varkappa^2 =0$;
\item the class of minimal surfaces of general type is characterized by $K^2-\varkappa^2>0$.
\end{itemize}

The class of minimal super-conformal surfaces in $\R^4$ is locally
equivalent to the class of holomorphic curves in $\C^2 \equiv
\R^4$ \cite{Ein}.

The inequality $K^2-\varkappa^2\geq 0$ for minimal surfaces also
follows from the inequality $K + |\varkappa| \leq ||H||^2$, which
holds for an arbitrary surface in $\R^4$ \cite{W}. Following
\cite{P-T&V}, a surface in $\R^4$ is called
\emph{Wintgen ideal surface}, if it satisfies the equality case of
the Wintgen's inequality identically. The Wintgen ideal surfaces
are characterized by circular ellipse of normal curvature
\cite{Guad-Rod}. In \cite{Chen2} B.-Y. Chen completely classified
Wintgen ideal surfaces in $\R^4$ with equal Gauss and normal
curvatures, i.e. Wintgen ideal surfaces satisfying $|K| =
|\varkappa|$ identically.

\vskip 2mm The surfaces with flat normal connection are
characterized in terms of the ellipse of normal curvature as
follows

\begin{prop}\label{P:flat-normal-ellipse}
Let $M^2$ be a surface in $\R^4$ free of flat points. Then $M^2$
is a surface with flat normal connection if and only if for each
point $p \in M^2$ the ellipse of normal curvature is a line
segment, which is not collinear with the mean curvature vector
field.
\end{prop}

\noindent {\it Proof:} The formula $S_{\varepsilon}=
\ds{\frac{\pi}{2}\, |\varkappa|}$ implies that $\varkappa = 0$ if
and only if the ellipse of curvature is a line segment.

Let $M^2$ be a surface with flat normal connection, i.e.
$\varkappa = 0$, $k \neq 0$. We assume that $M^2$ is parameterized
by principal parameters. Then from (2.2) it follows, that $\nu_1 =
\nu_2$. Further, equalities (2.1) imply that for each $v = \cos
\psi \, x + \sin \psi \, y$, we have $\sigma(v, v) = H + \sin 2
\psi (\lambda\,b + \mu\,l)$. So, when $v$ goes once around the
unit tangent circle, the vector $\sigma(v,v)$ goes twice along the
line segment collinear with $\lambda\,b + \mu\,l$ and centered at
$H$. The mean curvature vector field is $H = \nu_1\, b$.  Since $k
\neq 0$, then $\mu \neq 0$, and the line segment is not collinear
with $H$. \qed

\vskip 2mm We note that in the case $\lambda = 0$ the mean
curvature vector field $H$ is orthogonal to the line segment,
while in the case  $\lambda \neq 0$ the mean curvature vector
field $H$ is not orthogonal to the line segment. The length $d$ of
the line segment is
$$d = \sqrt{\lambda^2 + \mu^2} = \sqrt {H^2 - K}.$$
So, there arises a subclass of surfaces with flat normal
connection, characterized by the conditions:
$$ K = 0  \quad {\rm or} \quad d = \Vert H \Vert.$$

Proposition \ref{P:Flat normal-Lorentz circle} and Proposition
\ref{P:flat-normal-ellipse} give us the following

\begin{cor}\label{C:Lorentz circle-normal ellipse}
Let $M^2$ be a surface in $\R^4$ free of flat points.  Then the
tangent indicatrix $\chi$ is a rectangular hyperbola (a Lorentz
circle) if and only if the ellipse of normal curvature is a line
segment, which is not collinear with the mean curvature vector
field.
\end{cor}

\section{Fundamental theorem}

The basic theorem in the local differential geometry of surfaces
in $\R^3$ is the fundamental theorem of Bonnet. We proved a
theorem of Bonnet-type for minimal surfaces in $\R^4$ in terms of
their invariants \cite{GM2}. We recall that minimal surfaces are
characterized by the condition $\varkappa^2 - k = 0$ at any point.
In this section we consider surfaces in $\R^4$ free of minimal
points, i.e.
$$\varkappa^2 - k > 0 \qquad \textrm{at any point}.$$

We assume that $M^2$ is parameterized by principal lines and
denote the unit vector fields $\displaystyle{x=\frac{z_u}{\sqrt
E}, \; y=\frac{z_v}{\sqrt G}}$. If $\{x,y,b,l\}$ is the
geometrically determined moving frame field, we have the following
Frenet-type formulas:
$$\begin{array}{ll}
\vspace{2mm} \nabla'_xx=\quad \quad \quad \gamma_1\,y+\,\nu_1\,b;
& \qquad
\nabla'_xb=-\nu_1\,x-\lambda\,y\quad\quad \quad +\beta_1\,l;\\
\vspace{2mm} \nabla'_xy=-\gamma_1\,x\quad \quad \; + \; \lambda\,b
\; + \mu\,l;  & \qquad
\nabla'_yb=-\lambda\,x - \; \nu_2\,y\quad\quad \quad +\beta_2\,l;\\
\vspace{2mm} \nabla'_yx=\quad\quad \;-\gamma_2\,y \; + \lambda\,b
\; +\mu\,l;  & \qquad
\nabla'_xl= \quad \quad \quad \;-\mu\,y-\beta_1\,b;\\
\vspace{2mm} \nabla'_yy=\;\;\gamma_2\,x \quad\quad\quad+\nu_2\,b;
& \qquad \nabla'_yl=-\mu\,x \quad \quad \quad \;-\beta_2\,b,
\end{array}\leqno{(4.1)}$$
where  $\gamma_1 = - y(\ln \sqrt{E}), \,\, \gamma_2 = - x(\ln
\sqrt{G})$,
 $\nu_1, \nu_2, \lambda, \mu, \beta_1, \beta_2$  are geometric invariant functions.
Since $\varkappa^2 - k > 0$, then equalities (2.2) imply that $\mu
\neq 0$.

Using that $R'(x,y,x) = 0$, $R'(x,y,y) = 0$, $R'(x,y,b) = 0$ and
$R'(x,y,l) = 0$,  from (4.1) we get the following integrability
conditions:
$$\begin{array}{l}
\vspace{2mm}
2\mu\, \gamma_2 + \nu_1\,\beta_2 - \lambda\,\beta_1 = x(\mu);\\
\vspace{2mm}
2\mu\, \gamma_1 - \lambda\,\beta_2 + \nu_2\,\beta_1 = y(\mu);\\
\vspace{2mm}
\nu_1 \,\nu_2 - (\lambda^2 + \mu^2) = x(\gamma_2) + y(\gamma_1) - \left((\gamma_1)^2 + (\gamma_2)^2\right);\\
\vspace{2mm}
2\lambda\, \gamma_2 + \mu\,\beta_1 - (\nu_1 - \nu_2)\,\gamma_1 = x(\lambda) - y(\nu_1);\\
\vspace{2mm}
2\lambda\, \gamma_1 + \mu\,\beta_2 + (\nu_1 - \nu_2)\,\gamma_2 = - x(\nu_2) + y(\lambda);\\
\gamma_1\,\beta_1 - \gamma_2\,\beta_2 + (\nu_1 - \nu_2)\,\mu  = -
x(\beta_2) + y(\beta_1).
\end{array}\leqno{(4.2)}$$

Further we consider the general class of surfaces determined by
the condition
$$\mu_u \,\mu_v \neq 0.$$
The first two equalities of (4.2) imply that the condition
$\mu_u \,\mu_v \neq 0$ is equivalent to $$(2\mu\, \gamma_2 + \nu_1\,\beta_2
- \lambda\,\beta_1) (2\mu\, \gamma_1 - \lambda\,\beta_2 +
\nu_2\,\beta_1) \neq 0.$$
Then
$$\sqrt{E} = \displaystyle{\frac{\mu_u}{2\mu\, \gamma_2 + \nu_1\,\beta_2 - \lambda\,\beta_1}}, \quad
\sqrt{G} = \displaystyle{\frac{\mu_v}{2\mu\, \gamma_1 -
\lambda\,\beta_2 + \nu_2\,\beta_1}},$$
and the geometric invariants of the surface satisfy the inequalities
$$\displaystyle{\frac{\mu_u}{2\mu\, \gamma_2 + \nu_1\,\beta_2 - \lambda\,\beta_1}}>0, \quad
\displaystyle{\frac{\mu_v}{2\mu\, \gamma_1 - \lambda\,\beta_2 + \nu_2\,\beta_1}}>0. \leqno{(4.3)}$$

Furthermore, taking into account (4.2) we get that the invariants of any surface from the general class satisfy
the equalities
$$\begin{array}{l}
\vspace{2mm}
- \gamma_1 \sqrt{E} \sqrt{G} = (\sqrt{E})_v;\\
\vspace{2mm}
- \gamma_2 \sqrt{E} \sqrt{G} = (\sqrt{G})_u;\\
\vspace{2mm} \nu_1 \,\nu_2 - (\lambda^2 + \mu^2) =
\displaystyle{\frac{1}{\sqrt{E}}\,(\gamma_2)_u
+ \frac{1}{\sqrt{G}}\,(\gamma_1)_v - \left((\gamma_1)^2 + (\gamma_2)^2\right)};\\
\vspace{2mm} 2\lambda\, \gamma_2 + \mu\,\beta_1 - (\nu_1 -
\nu_2)\,\gamma_1 =
\displaystyle{\frac{1}{\sqrt{E}}\,\lambda_u - \frac{1}{\sqrt{G}}\,(\nu_1)_v};\\
\vspace{2mm} 2\lambda\, \gamma_1 + \mu\,\beta_2 + (\nu_1 -
\nu_2)\,\gamma_2 =
\displaystyle{ - \frac{1}{\sqrt{E}}\,(\nu_2)_u + \frac{1}{\sqrt{G}}\,\lambda_v};\\
\gamma_1\,\beta_1 - \gamma_2\,\beta_2 + (\nu_1 - \nu_2)\,\mu  =
\displaystyle{ - \frac{1}{\sqrt{E}}\,(\beta_2)_u +
\frac{1}{\sqrt{G}}\,(\beta_1)_v}.
\end{array}\leqno{(4.4)}$$

\vskip 3mm We shall prove the following Bonet-type fundamental
theorem

\begin{thm}\label{T:Main Theorem}
Let $\gamma_1, \, \gamma_2, \, \nu_1,\, \nu_2, \, \lambda, \, \mu,
\, \beta_1, \beta_2$ be smooth functions, defined in a domain
$\mathcal{D}, \,\, \mathcal{D} \subset {\R}^2$, and satisfying
inequalities {\rm (4.3)} and  equalities {\rm (4.4)}.

Let $x_0, \, y_0, \, b_0,\, l_0$ be an
orthonormal frame at a point $p_0 \in \R^4$. Then there exist a
subdomain ${\mathcal{D}}_0 \subset \mathcal{D}$ and a unique
surface $M^2: z = z(u,v), \,\, (u,v) \in {\mathcal{D}}_0$, passing
through $p_0$, such that $\gamma_1, \, \gamma_2, \, \nu_1,\,
\nu_2, \, \lambda, \, \mu, \, \beta_1, \beta_2$ are the geometric
functions of $M^2$ and $x_0, \, y_0, \, b_0,\, l_0$ is the
geometric frame of $M^2$ at the point $p_0$.
\end{thm}
\vskip 2mm \noindent \emph{Proof:} We consider the following
system of partial differential equations for the unknown vector
functions $x = x(u,v), \, y = y(u,v), \,b = b(u,v), \,l = l(u,v)$
in $\R^4$:
$$\begin{array}{ll}
\vspace{2mm} x_u = \sqrt{E}\, \gamma_1\, y + \sqrt{E}\, \nu_1\, b
& \qquad x_v =
- \sqrt{G}\, \gamma_2\, y + \sqrt{G}\, \lambda\, b + \sqrt{G}\, \mu\, l\\
\vspace{2mm} y_u = - \sqrt{E}\, \gamma_1\, x + \sqrt{E}\,
\lambda\, b
+ \sqrt{E}\, \mu\, l  & \qquad y_v = \sqrt{G}\, \gamma_2\, x + \sqrt{G}\, \nu_2\, b \\
\vspace{2mm} b_u = - \sqrt{E}\, \nu_1\, x - \sqrt{E}\, \lambda\, y
+ \sqrt{E}\, \beta_1\, l  & \qquad b_v =
- \sqrt{G}\, \lambda\, x - \sqrt{G}\, \nu_2\, y + \sqrt{G}\, \beta_2\, l \\
\vspace{2mm} l_u = - \sqrt{E}\, \mu\, y - \sqrt{E}\, \beta_1\, b &
\qquad l_v = - \sqrt{G}\, \mu\, x - \sqrt{G}\, \beta_2\, b
\end{array}\leqno{(4.5)}$$
We denote
$$Z =
\left(%
\begin{array}{c}
  x \\
  y \\
  b \\
  l \\
\end{array}%
\right); \quad
A = \sqrt{E} \left(%
\begin{array}{cccc}
  0 & \gamma_1 & \nu_1 & 0 \\
  -\gamma_1 & 0 & \lambda & \mu \\
  -\nu_1 & -\lambda & 0 & \beta_1 \\
  0 & -\mu & -\beta_1 & 0 \\
\end{array}%
\right); \quad B = \sqrt{G}
\left(%
\begin{array}{cccc}
  0 & -\gamma_2 & \lambda & \mu \\
  \gamma_2 & 0 & \nu_2 & 0 \\
  -\lambda & -\nu_2 & 0 & \beta_2 \\
  -\mu & 0 & -\beta_2 & 0 \\
\end{array}%
\right).$$ Then system (4.5) can be rewritten in the form:
$$\begin{array}{l}
\vspace{2mm}
Z_u = A\,Z,\\
\vspace{2mm} Z_v = B\,Z.
\end{array}\leqno{(4.6)}$$
The integrability conditions of (4.6) are
$$Z_{uv} = Z_{vu},$$
i.e.
$$\displaystyle{\frac{\partial a_i^k}{\partial v} - \frac{\partial b_i^k}{\partial u}
+ \sum_{j=1}^{4}(a_i^j\,b_j^k - b_i^j\,a_j^k) = 0, \quad i,k = 1,
\dots, 4,} \leqno(4.7)$$ where $a_i^j$ and $b_i^j$ are the
elements of the matrices $A$ and $B$. Using (4.4) we obtain that
the equalities (4.7) are fulfilled. Hence, there exists a subset
$\mathcal{D}_1 \subset \mathcal{D}$ and unique vector functions $x
= x(u,v), \, y = y(u,v), \,b = b(u,v), \,l = l(u,v), \,\, (u,v)
\in \mathcal{D}_1$, which satisfy system (4.5) and the conditions
$$x(u_0,v_0) = x_0, \quad y(u_0,v_0) = y_0, \quad b(u_0,v_0) = b_0, \quad l(u_0,v_0) = l_0.$$

We shall prove that $x(u,v), \, y(u,v), \,b(u,v), \,l(u,v)$ form
an orthonormal frame for each $(u,v) \in \mathcal{D}_1$. Let us
consider the following functions:
$$\begin{array}{lll}
\vspace{2mm}
  \varphi_1 = x^2 - 1; & \qquad \varphi_5 = x\,y; & \qquad \varphi_8 = y\,b; \\
\vspace{2mm}
  \varphi_2 = y^2 - 1; & \qquad \varphi_6 = x\,b; & \qquad \varphi_9 = y\,l; \\
\vspace{2mm}
  \varphi_3 = b^2 - 1; & \qquad \varphi_7 = x\,l; & \qquad \varphi_{10} = b\,l; \\
\vspace{2mm}
  \varphi_4 = l^2 - 1; &   &  \\
\end{array}$$
defined for each $(u,v) \in \mathcal{D}_1$. Using that $x(u,v), \,
y(u,v), \,b(u,v), \,l(u,v)$ satisfy (4.5), we obtain  the system
$$\begin{array}{lll}
\vspace{2mm}
\displaystyle{\frac{\partial \varphi_i}{\partial u} = \alpha_i^j \, \varphi_j},\\
\vspace{2mm} \displaystyle{\frac{\partial \varphi_i}{\partial v} =
\beta_i^j \, \varphi_j};
\end{array} \qquad i = 1, \dots, 10, \leqno{(4.8)}$$
where $\alpha_i^j, \beta_i^j, \,\, i,j = 1, \dots, 10$ are
functions of $(u,v) \in \mathcal{D}_1$. System (4.8) is a linear
system of partial differential equations for the functions
$\varphi_i(u,v), \,\,i = 1, \dots, 10, \,\,(u,v) \in
\mathcal{D}_1$, satisfying $\varphi_i(u_0,v_0) = 0, \,\,i = 1,
\dots, 10$. Hence, $\varphi_i(u,v) = 0, \,\,i = 1, \dots, 10$ for
each $(u,v) \in \mathcal{D}_1$. Consequently, the vector functions
$x(u,v), \, y(u,v), \,b(u,v), \,l(u,v)$ form an orthonormal frame
field for each $(u,v) \in \mathcal{D}_1$.

Now, let us consider the  system
$$\begin{array}{lll}
\vspace{2mm}
z_u = \sqrt{E}\, x\\
\vspace{2mm} z_v = \sqrt{G}\, y
\end{array}\leqno{(4.9)}$$
of partial differential equations for the vector function
$z(u,v)$. Using (4.4) and (4.5) we get that the integrability
conditions $z_{uv} = z_{vu}$ of system (4.9)
 are fulfilled. Hence,  there exists a subset $\mathcal{D}_0 \subset \mathcal{D}_1$ and
a unique vector function $z = z(u,v)$, defined for $(u,v) \in
\mathcal{D}_0$ and satisfying $z(u_0, v_0) = p_0$.

Consequently, the surface $M^2: z = z(u,v), \,\, (u,v) \in
\mathcal{D}_0$ satisfies the assertion of the theorem. \qed

\section{Examples of surfaces with flat normal connection} \label{S:Examples}

In this section we construct a family of surfaces with flat normal
connection lying on a standard rotational hypersurface in $\R^4$ .

Let $\{e_1, e_2, e_3, e_4\}$ be the standard orthonormal frame in
$\R^4$, and $S^2(1)$ be a 2-dimensional sphere in $\R^3 = \span
\{e_1, e_2, e_3\}$, centered at the origin $O$.
Let $f = f(u), \,\, g = g(u)$ be smooth functions, defined in an
interval $I \subset \R$, such that $\dot{f}^2(u) + \dot{g}^2(u) =
1, \,\, u \in I$.
The standard rotational hypersurface $M^3$ in
$\R^4$, obtained by the rotation of the meridian curve $m: u
\rightarrow (f(u), g(u))$ around the $Oe_4$-axis,  is parameterized as follows:
$$M^3: Z(u,w^1,w^2) = f(u)\,l(w^1,w^2) + g(u) \,e_4,$$
where $l(w^1,w^2)$ is the unit radius-vector of $S^2(1)$  in $\R^3$.

We consider a
smooth curve $c: l = l(v) = l(w^1(v),w^2(v)), \, v \in J, \,\, J \subset \R$  on
$S^2(1)$, parameterized by the arc-length, i.e. $l'^2(v) = 1$. We
denote $t = l'$ and consider the moving frame field $\span \{t(v),
n(v), l(v)\}$ of the curve $c$ on $S^2(1)$. With respect to this
orthonormal frame field the following Frenet formulas hold good:
$$\begin{array}{l}
\vspace{2mm}
l' = t;\\
\vspace{2mm}
t' = \kappa \,n - l;\\
\vspace{2mm} n' = - \kappa \,t,
\end{array} \leqno{(5.1)}$$
where $\kappa$ is the spherical curvature of $c$.

 Now we construct a surface $M^2$ in $\R^4$ in
the following way:
$$M^2: z(u,v) = f(u) \, l(v) + g(u)\, e_4, \quad u \in I, \, v \in J. \leqno{(5.2)}$$
The surface $M^2$ lies on the rotational hypersurface $M^3$ in
$\R^4$. Since $M^2$ consists of meridians of $M^3$, we call $M^2$ a
\textit{meridian surface}.

The tangent space of $M^2$ is spanned by the vector fields:
$$\begin{array}{l}
\vspace{2mm}
z_u = \dot{f} \,l + \dot{g}\,e_4;\\
\vspace{2mm} z_v = f\,t,
\end{array}$$
and hence, the coefficients of the first fundamental form of $M^2$
are $E = 1; \,\, F = 0; \,\, G = f^2(u)$. Taking into account
(5.1), we calculate the second partial derivatives of $z(u,v)$:
$$\begin{array}{l}
\vspace{2mm}
z_{uu} = \ddot{f} \,l + \ddot{g}\,e_4;\\
\vspace{2mm}
z_{uv} = \dot{f}\,t;\\
\vspace{2mm} z_{vv} = f \kappa \,n - f\,l.
\end{array}$$
Let us denote $x = z_u,\,\, y = \ds{\frac{z_v}{f} = t}$ and
consider the following orthonormal normal frame field of $M^2$:
$$n_1 = n(v); \qquad n_2 = - \dot{g}(u)\,l(v) + \dot{f}(u) \, e_4.$$
Thus we obtain a positive orthonormal frame field $\{x,y, n_1,
n_2\}$ of $M^2$. If we denote by $\kappa_m$ the curvature of the
meridian curve $m$, i.e. $\kappa_m (u)= \dot{f}(u) \ddot{g}(u) -
\dot{g}(u) \ddot{f}(u) = \ds{\frac{- \ddot{f}(u)}{\sqrt{1 -
\dot{f}^2(u)}}}$, then we get the following derivative formulas of
$M^2$:
$$\begin{array}{ll}
\vspace{2mm} \nabla'_xx = \qquad \qquad \qquad \qquad
\kappa_m\,n_2; & \qquad
\nabla'_x n_1 = 0;\\
\vspace{2mm} \nabla'_xy = 0;  & \qquad
\nabla'_y n_1 = \ds{\quad \quad \quad - \frac{\kappa}{f}\,y};\\
\vspace{2mm} \nabla'_yx = \quad\quad
\quad\ds{\frac{\dot{f}}{f}}\,y;  & \qquad
\nabla'_x n_2 = - \kappa_m \,x;\\
\vspace{2mm} \nabla'_yy = \ds{- \frac{\dot{f}}{f}\,x \quad\quad +
\frac{\kappa}{f}\,n_1 + \frac{\dot{g}}{f} \, n_2}; & \qquad
\nabla'_y n_2 = \ds{ \quad \quad \quad - \frac{\dot{g}}{f}\,y}.
\end{array} \leqno{(5.3)}$$

The coefficients of the second fundamental form of $M^2$ are $L =
N = 0, \,\, M = - \kappa_m(u) \, \kappa(v)$. Taking into account
(5.3), we find the invariants $k$, $\varkappa$, $K$:
$$k = - \frac{\kappa_m^2(u) \, \kappa^2(v)}{f^2(u)}; \qquad \varkappa = 0;
\qquad K = \frac{\kappa_m (u)\, \dot{g}(u)}{f(u)}. \leqno{(5.4)}$$

The equality $\varkappa = 0$ implies that $M^2$ is a surface with
flat normal connection.

The mean curvature vector field $H$ is given by
$$H = \frac{\kappa}{2f}\, n_1 + \frac{\dot{g} + f \kappa_m}{2f} \, n_2. \leqno{(5.5)}$$

There are three main classes of meridian surfaces: \vskip 2mm I.
$\kappa = 0$, i.e. the curve $c$ is a great circle on $S^2(1)$. In
this case $n_1 = const$, and $M^2$ is a planar surface lying in
the constant 3-dimensional space spanned by $\{x, y, n_2\}$.
Particularly, if in addition $\kappa_m = 0$, i.e. the meridian
curve lies on a straight line, then $M^2$ is a developable surface
in the 3-dimensional space $\span \{x, y, n_2\}$.

\vskip 2mm II. $\kappa_m = 0$, i.e. the meridian curve is part of
a straight line. In such case $k = \varkappa = K = 0$, and $M^2$
is a developable ruled surface. If in addition $\kappa = const$,
i.e. $c$ is a circle on $S^2(1)$, then $M^2$ is a developable
ruled surface in a 3-dimensional space. If $\kappa \neq const$,
i.e. $c$ is not a circle on $S^2(1)$, then $M^2$ is a developable
ruled surface in $\R^4$.

\vskip 2mm III. $\kappa_m \, \kappa \neq 0$, i.e. $c$ is not a
great circle on $S^2(1)$, and $m$ is not a straight line. In this
general case the invariant function $k<0$, which implies that
there exist two systems of asymptotic lines on $M^2$. The
parametric lines of $M^2$ given by (5.2) are orthogonal and
asymptotic.

\vskip 2mm Let $M^2$ be a meridian surface of the general class.
Now we are going to find the meridian surfaces with:
\begin{itemize}
\item
constant Gauss curvature $K$;
\item
constant mean curvature;
\item
constant invariant function $k$.
\end{itemize}

\begin{prop}
Let $M^2$ be a meridian surface in $\R^4$ from the general class.
Then $M^2$ has constant non-zero Gauss curvature $K$ if and only
if the meridian $m$ is given by
$$\begin{array}{ll}
\vspace{2mm}
f(u) = \alpha \cos \sqrt{K} u + \beta \sin \sqrt{K} u, & \quad K >0;\\
\vspace{2mm} f(u) = \alpha \cosh \sqrt{-K} u + \beta \sinh
\sqrt{-K} u, & \quad K <0,
\end{array}$$
where $\alpha$ and $\beta$ are constants.
\end{prop}

\noindent {\it Proof:} Using (5.4) and $\dot f^2+\dot g^2= 1$, we
obtain that $M^2$ has constant Gauss curvature $K \neq 0$ if and
only if the meridian $m$ satisfies the following differential
equation
$$\ddot{f}(u) + K f(u) = 0.$$
The general solution of the above equation is given by
$$\begin{array}{ll}
\vspace{2mm}
f(u) = \alpha \cos \sqrt{K} u + \beta \sin \sqrt{K} u, & \quad \textrm{in case} \quad K >0;\\
\vspace{2mm} f(u) = \alpha \cosh \sqrt{-K} u + \beta \sinh
\sqrt{-K} u, & \quad \textrm{in case}  \quad K <0,
\end{array}$$
where $\alpha$ and $\beta$ are constants. The function $g(u)$ is
determined by $\dot{g}(u) = \sqrt{1 - \dot{f}^2(u)}$.

\qed

\vskip 3mm The equality (5.5) implies that the mean curvature of
$M^2$ is given by
$$|| H || = \ds{\sqrt{\frac{\kappa^2(v) + \left(\dot{g}(u)
+ f(u) \kappa_m(u)\right)^2}{4f^2(u)}}}. \leqno{(5.6)}$$

The  meridian surfaces with constant mean curvature (CMC meridian
surfaces) are described in

\begin{prop}
Let $M^2$ be a meridian surface in $\R^4$ from the general class.
Then $M^2$ has constant mean curvature $|| H || = a = const$, $a
\neq 0$ if and only if the curve $c$  on $S^2(1)$ is a circle with
constant spherical curvature $\kappa = const = b, \; b \neq 0$,
and the meridian $m$ is determined by the following differential
equation:
$$\left( 1 - \dot{f}^2 - f \ddot{f}\right)^2 = (1 - \dot{f}^2) (4 a^2 f^2 - b^2).$$
\end{prop}

\noindent {\it Proof:} From (5.6) it follows that $||H|| = a$ if
and only if
$$\kappa^2(v) = 4 a^2 f^2(u) - (\dot{g}(u) + f(u)  \kappa_m(u))^2,$$
which implies
$$\begin{array}{l}
\vspace{2mm}
\kappa = const = b, \; b \neq 0;\\
\vspace{2mm} 4 a^2 f^2(u) - (\dot{g}(u) + f(u)  \kappa_m(u))^2 =
b^2.
\end{array} \leqno{(5.7)}$$
The first equality of (5.7) implies that the spherical curve $c$
has constant spherical curvature $\kappa = b$,
 i.e. $c$ is a circle.
Using that $\dot{f}^2 + \dot{g}^2 = 1$, and $\kappa_m = \dot{f}
\ddot{g} - \dot{g} \ddot{f}$ we calculate that $\dot{g} + f
\kappa_m = \ds{\frac{1- \dot{f}^2 - f \ddot{f}}{\sqrt{1 -
\dot{f}^2}}}.$ Hence, the second equality of (5.7) gives the
following differential equation for the meridian $m$:
$$\left( 1 - \dot{f}^2 - f \ddot{f}\right)^2 = (1 - \dot{f}^2) (4 a^2 f^2 - b^2). \leqno{(5.8)}$$

Further, if we set $\dot{f} = y(f)$ in equation (5.8), we obtain
that the function $y = y(t)$ is a solution of the following
differential equation
$$1 - y^2 - \frac{t}{2}(y^2)' = \sqrt{1 - y^2} \sqrt{4 a^2 t^2 - b^2}.$$
The general solution of the above equation is given by
$$y(t) = \sqrt{1 - \frac{1}{t^2}\left( C + \frac{t}{2} \sqrt{4 a^2 t^2 - b^2} -
\frac{b^2}{4a} \ln |2at + \sqrt{4 a^2 t^2 - b^2}| \right)^2};
\qquad C = const. \leqno{(5.9)}$$ The function $f(u)$ is
determined by $\dot{f} = y(f)$ and (5.9). The function $g(u)$ is
defined by $\dot{g}(u) = \sqrt{1 - \dot{f}^2(u)}$. \qed

\vskip 3mm At the end of this section we shall find the meridian
surfaces with constant invariant $k$.

\begin{prop}
Let $M^2$ be a meridian surface in $\R^4$ from the general class.
Then $M^2$ has a constant invariant $k = const = - a^2, \; a\neq
0$ if and only if the curve $c$  on $S^2(1)$ is a circle with
spherical curvature $\kappa = const = b, \; b \neq 0$, and the
meridian $m$ is determined by the following differential equation:
$$\ddot{f}(u) = \mp \frac{a}{b}\,f(u) \sqrt{1 - \dot{f}^2(u)}.$$
\end{prop}

\noindent {\it Proof:} Using (5.4) we obtain that $k = const = -
a^2, \; a \neq 0$ if and only if $\kappa^2(v) \kappa_m^2(u) = a^2
f^2(u)$. Hence,
$$\kappa(v) = \pm \, a \, \frac{f(u)}{\kappa_m(u)}.$$
The last equality implies
$$\begin{array}{l}
\vspace{2mm}
\kappa = const = b, \; b \neq 0;\\
\vspace{2mm}
 \pm \, a \, \ds{\frac{f(u)}{\kappa_m(u)}} = b.
\end{array} \leqno{(5.10)}$$
The first equality of (5.10) implies that the spherical curve $c$
has constant spherical curvature $\kappa = b$,  i.e. $c$ is a
circle. The second equality of (5.10) gives the following
differential equation for the function $f(u)$:
$$\frac{\ddot{f}(u)}{\sqrt{1 - \dot{f}^2(u)}} = \mp \frac{a}{b}\,f(u).  \leqno{(5.11)}$$
Again setting $\dot{f} = y(f)$ in equation (5.11), we obtain that
the function $y = y(t)$ is a solution of the following
differential equation

$$\frac{y y'}{\sqrt{1 - y^2}} = \mp \frac{a}{b}\,t.$$
The general solution of the above equation is given by
$$y(t) = \sqrt{1 - \left( C \pm \frac{a}{b}\, \frac{t^2}{2} \right)^2}; \qquad C = const.  \leqno{(5.12)}$$
The function $f(u)$ is determined by $\dot{f} = y(f)$ and (5.12).
The function $g(u)$ is defined by $\dot{g}(u) = \sqrt{1 -
\dot{f}^2(u)}$. \qed \vskip 2mm

 \vskip 3mm \textbf{Acknowledgements:}
The authors would like to
express their thanks to the referees for their valuable comments and
suggestions.
  The second author is
partially supported by "L. Karavelov" Civil Engineering Higher
School, Sofia, Bulgaria under Contract No 10/2010.


\vskip 10mm
\begin{thebibliography}{99}

\bibitem{Asp}
Asperti A.,  Some generic properties of Riemannian immersions,
Bol. Soc. Brasil. Mat., 1980,  11, no. 2, 191--216.


\bibitem{BFLPP}
Burstall F., Ferus D.,  Leschke K.,  Pedit F., Pinkall  U.,
Conformal geometry of surfaces in the 4-sphere and quaternions,
Lecture Notes in Mathematics, vol. 1772, Springer-Verlag,
2002.

\bibitem{Chen1}
Chen B.-Y.,  Geometry of submanifolds, Marcel Dekker, Inc.,
New York, 1973.

\bibitem{Chen2}
Chen B.-Y.,  Classification of Wintgen ideal surfaces in Euclidean
4-space with equal Gauss and normal curvatures,  Ann. Glob.
Anal. Geom., 2010, 38, 145--160.


\bibitem{DT}
Dajczer M., Tojeiro R., All superconformal surfaces in ${\R^4}$
in terms of minimal surfaces,  Math. Z., 2009, 261,
no. 4, 869--890.

\bibitem{Ein}
Eisenhart L., Minimal surfaces in Euclidean four-space,  Amer.
J. Math., 1912, 34, 215--236.

\bibitem{GM1}
Ganchev G., Milousheva V.,  On the theory of surfaces in the
four-dimensional Euclidean space,   Kodai Math. J.,
2008,   31, 183--198.

\bibitem{GM2}
Ganchev G., Milousheva V.,  Minimal surfaces in the
four-dimensional Euclidean space, preprint available at http://arxiv.org/abs/0806.3334v1

\bibitem{GM3}
Ganchev G., Milousheva V., Invariants of lines on surfaces in
$\R^4$,  C. R. Acad. Bulg. Sci., 2010,  63, no. 6,
835--842.

\bibitem{Gar-Soto}
Garcia R., Sotomayor J., Lines of axial curvatures on surfaces
immersed in $\R^4$,  Differential Geom. Appl.,
2000,  12, 253--269.

\bibitem{GVV1}
Gheysens L.,  Verheyen P.,  Verstraelen L., Sur les surfaces
$\mathcal{A}$ ou les surfaces de Chen,   C. R. Acad. Sci.
Paris, S\'{e}r. I, 1981, 292, 913--916.


\bibitem{GVV2}
Gheysens L.,  Verheyen P.,  Verstraelen L., Characterization
and examples of Chen submanifolds,  J. Geom.,
1983,  20, 47--62.


\bibitem{Guad-Rod}
Guadalupe I., Rodriguez L., Normal curvature of surfaces in
space forms, Pacific J. Math.,  1983, 106, no. 1,
95--102.


\bibitem{Kom}
Kommerell K., Riemannsche Fl\"{a}chen im ebenen Raum von vier
Dimensionen, Math. Ann., 1905, 60, 546--596.

\bibitem{Little}
Little J., On singularities of submanifolds of higher dimensional
Euclidean spaces, Ann. Mat. Pura Appl., IV. Ser, 1969, 83,
261--335.


\bibitem{Mello}
Mello L., Orthogonal asymptotic lines on surfaces immersed in
$\mathbb R^4$,  Rocky Mountain J. Math.,
2009, 39, 1597--1612.

\bibitem{MW1}
Moore C.,  Wilson E., A general theory of surfaces,  J.
Nat. Acad. Proc., 1916, 2, 273--278.

\bibitem{MW2}
Moore C.,  Wilson E.,  Differential geometry of two-dimensional
surfaces in hyperspaces,   Proc. Acad. Arts Sci.,
1916, 52, 267--368.

\bibitem{P-T&V}
Petrovi\'{c}-Torga\v{s}ev M.,  Verstraelen L.,  On Deszcz
symmetries of Wintgen ideal submanifolds, Arch. Math. (Brno),
2008, 44, 57--76.

\bibitem{Scho-Str}
Schouten J.,  Struik D., Einf\"{u}hrung in die neueren
Methoden der Differentialgeometrie II, Batavia, 1938.

\bibitem{Spivak}
 Spivak M., Introduction to Comprehensive Differential
Geometry, vol. I, V, Publish or Perish, Berkeley, 1999.


\bibitem{W}
Wintgen P., S\"{u}r l'inegalit\'{e} de Chen-Willmore, C. R.
Acad. Sc. Paris, S\'{e}r. A, 1979, 288, 993-995.

\bibitem{Wong}
 Wong Y.-C.,  A new curvature theory for surfaces in a Euclidean
4-space,  Comm. Math. Helv., 1952, 26, 152--170.


\end{thebibliography}
\end{document}